\documentclass[12pt,a4paper]{article}
\usepackage[T1]{fontenc}
\usepackage[british]{babel}
\usepackage[utf8]{inputenc}
\usepackage{amsmath,amssymb,amsthm}
\usepackage{hyperref}
\newtheorem{thm}{Theorem}[section]

\newtheorem{lem}[thm]{Lemma}

\usepackage{chngcntr}     

\newtheorem{theorem}{Theorem}

\counterwithout{theorem}{section}  

\DeclareMathOperator{\Irr}{Irr}

\DeclareMathOperator{\Ker}{ker}
\DeclareMathOperator{\PSL}{PSL}

\DeclareMathOperator{\N}{N}
\DeclareMathOperator{\M}{M}
\DeclareMathOperator{\SL}{SL}
\DeclareMathOperator{\Core}{Core}
\DeclareMathOperator{\Oo}{O}
\DeclareMathOperator{\Cent}{C}
\title{On permutation characters of finite groups}
\author{Jiakuan Lu\thanks{School of Mathematics and Statistics, Guangxi Normal University,
		Guilin 541006, Guangxi,  P. R. China.}~~and  Hangyang Meng\thanks{Email:
		\texttt{hymeng2009@shu.edu.cn}. Department of Mathematics, Shanghai University, and
	Newtouch Center for Mathematics of Shanghai University,
	Shanghai 200444, P. R. China.}
	}
\begin{document}
\date{Dedicated to the memory of I. Martin Isaacs}
\maketitle
\begin{abstract}
Let $G$ be a finite group and \( M \) be a maximal subgroup of \( G \).  We say that  every irreducible constituent \( \chi \)  of \( (1_M)^G \) is a \( \mathcal{P} \)-character of \( G \) with respect to \( M \). In this paper, we prove that  $G$ is solvable if and only if all $\mathcal{P}$-characters of $G$ are monomial. This result is a generalization of Taketa's theorem on $\M$-groups and solves a question posed
by Qian and Yang.


\end{abstract}

\section{Introduction}

\maketitle
In this paper, $G$ will always denote a finite group. Let \( \text{Irr}(G) \) be the set of all irreducible (complex) characters of \( G \). Recall that a character $\chi$ of $G$ is monomial if it is induced from a linear character of some (not necessarily proper) subgroup, and we say that \( G \) is an \( \M \)-group if  every irreducible character of $G$ is monomial. It is well-known that \(\M \)-groups are solvable due to Taketa (see \cite[Corollary 5.13]{Isaacs1}).

 Monomial characters and $\M$-groups are probably one of the research topics that Prof. I. M. Isaacs was interested in during his lifetime. For example, in one of his early papers \cite{Isaacs84}, he obtained some generalizations of Taketa's theorem on the solvability of $\M$-groups, and Part 3 of his book \cite{Isaacs2} is primarily concerned with some of the deeper properties of $\M$-groups and monomial characters.


 
In this paper, for a character $\psi$ of a group $G$, we denote by $\Irr(\psi)=\{\chi \in \Irr(G) \mid [\chi,\psi] \neq 0\}$ the set of all irreducible constituents of $\psi$. Followed by Qian and Yang's paper~\cite{Qian}, we say that \( \chi \) is a \( \mathcal{P} \)-character of \( G \) with respect to the maximal subgroup \( M \) if \( \chi \in \Irr((1_M)^G)\).  
Based on the solvability of $\M$-groups, Qian and Yang also posed the following question \cite[Question~3.9]{Qian}: "Assume that all $\mathcal{P}$-characters of $G$ is monomial. Is $G$ solvable?"

Our first result will give an affirmative answer. In fact, we prove: 
\begin{theorem}\label{thma} 
Let $G$ be a group. Then the following statements are equivalent:
\begin{itemize}
\item[(1)] $G$ is solvable;
\item[(2)] All $\mathcal{P}$-characters of $G$ are monomial;
\item[(3)] $(1_M)^G$ has a non-principle, monomial, irreducible constituent for each non-normal maximal subgroup $M$ of $G$.
\end{itemize}
\end{theorem}

Gagola and Lewis in \cite{Gagola} proved that a group $G$ is nilpotent if and only if $\chi(1)^2$ divides $|G : \operatorname{ker}\chi|$ for each $\chi \in \Irr(G)$. In the following, we generalize this result to $\mathcal{P}$-characters in solvable groups. Here \( \text{Irr}_{\mathcal{P}}(G) \) denotes the set of all \( \mathcal{P} \)-characters of \( G \). 
\begin{theorem}\label{thmgl} 
Let $G$ be a solvable group. Then $G$ is nilpotent if and only if $\chi(1)^2$ divides $|G : \operatorname{ker}\chi|$ for each $\chi \in \operatorname{Irr}_{\mathcal{P}}(G)$.
\end{theorem}

We do not know that whether $G$ is solvable if $\chi(1)^2$ divides $|G : \operatorname{ker}\chi|$ for every character $\chi \in \operatorname{Irr}_{\mathcal{P}}(G)$. It is worth point out \cite[Corollary 2.6]{Qian} that a group $G$ is nilpotent if and only if each $\mathcal{P}$-character is linear.

  
  
  
 
 
For a set $\pi$ of some primes,  a character $\chi$ of $G$ is called a  \( \mathcal{P}_\pi \)-character if $\chi \in \Irr((1_M)^G)$ for some maximal subgroup $M$ of $G$ such that \( |G:M| \) is a \( \pi \)-number. The following theorem is a $\pi$-variation of  \cite[Proposition 2.5]{Qian}.
 
\begin{theorem}\label{thmb}
  Let $G$ be a $\pi$-solvable group. Then all ${\mathcal{P}_\pi}$-character degrees of $G$ are 
 $\pi$-number if and only if $G$ has a normal $\pi$-complement.
\end{theorem} 

Denote by $\Irr_{\mathcal{P}_{\pi}}(G)$ the set of all  \( \mathcal{P}_\pi \)-characters of $G$. Clearly
$$\Irr_{\mathcal{P}_{\pi}}(G) \subseteq \Irr_{\mathcal{P}}(G) \subseteq \Irr(G).$$
The next result characterizes \( \pi \)-separable groups satisfying $\Irr_{\mathcal{P}_{\pi}}(G)=\Irr_{\mathcal{P}}(G)$.
\begin{theorem}\label{thmc}
Let $\pi$ be a set of some primes and let \( G \) be a \( \pi \)-separable group. Then $\Irr_{\mathcal{P}_{\pi}}(G)=\Irr_{\mathcal{P}}(G)$ if and only if \( G \) is a \( \pi \)-group.
\end{theorem} 

Note that the "$\pi$-separability" of $G$ in Theorem~\ref{thmc} can not be removed. For example, $G=\PSL(2,7)$, the simple group of order 168, has maximal subgroups of indices $7$ and $8$. Then every $\mathcal{P}$-character of $G$ is a $\mathcal{P}_{\{2,7\}}$-character but $G$ is not a $\{2,7\}$-group.  

Note that $\Irr_{\mathcal{P}}(G) \subsetneq \Irr(G)$ in general. For example, $\SL(2,3)$ has irreducible characters of degree $2$ which are not $\mathcal{P}$-characters. Hence it will be interesting to investigate the class of groups satisfying $\Irr_{\mathcal{P}}(G)=\Irr(G)$.
As pointed out by Prof. G. Navarro [private
communication], such groups are not necessarily solvable, and some simple groups, such as ${\rm PSL}(2,27), M_{23},M_{24}$,  all satisfy \(\text{Irr}_{\mathcal{P}}(G) ={\rm Irr}(G) \).







\section{Proofs }
The first key lemma about the kernels of $\mathcal{P}$-characters is due to Qian and Yang~\cite{Qian}. Recall that, for a subgroup $M$ of a group $G$, we denote by $\Core_G(M)=\bigcap_{x \in G}M^x$ the core of $M$ in $G$, which is the largest normal subgroup of $G$ contained in $M$.
\begin{lem}{\rm \cite[Proposition 2.1]{Qian}}\label{lem-core}
Let \( \chi \) be a non-principal \( \mathcal{P} \)-character of \( G \) with respect to a maximal subgroup \( M \). Then \( \ker \chi ={\rm Core}_G(M). \) Moreover, $$\Phi(G)= \bigcap_{\chi \in \Irr_{\mathcal{P}}(G)} \ker \chi.$$
 \end{lem}
 


The following lemma is a slight generalization of partial results in \cite[Proposition~2.3]{Qian}. We give a proof here for completeness.
\begin{lem}\label{lem-monomial}
Let $G=NM$ be the semidirect product of subgroups $N$ and $M$ with $N \unlhd G$. Suppose that $\chi \in \Irr((1_M)^G)$ with $[\chi_N, \lambda] \neq 0$ for some linear character $\lambda \in \Irr(N)$. Then 
\begin{itemize}
\item[(1)]  $\chi(1)=|G:\Cent_G(\lambda)|$, where $\Cent_G(\lambda)=\{g \in G \mid \lambda^g=\lambda\}$; and 
\item[(2)]  $\chi=\psi^G$ for some $\psi \in \Irr(\Cent_G(\lambda))$ with $\psi_N=\lambda$. In particular, $\chi$ is monomial.
\end{itemize}
\end{lem}
\begin{proof}
As $M \cap N=1$, $((1_M)^G)_N=(1_{M \cap N})^N$ is a regular character of $N$, which implies that
$$((1_M)^G)_N=\sum_{\eta \in \Irr(N)} \eta(1)\eta.$$
Write $I=\Cent_G(\lambda)$. By Clifford's Theorem, $\chi_N=e\sum_{i=1}^m\lambda^{g_i}$, where $m=|G:I|$ and $\{g_1=1,g_2,\cdots,g_m\}$ is a transversal of $I$ in $G$. Note that, as $\lambda$ is linear,
$$1 \leq e=[\chi_N,\lambda] \leq [((1_M)^G)_N, \lambda]=[\sum_{\eta \in \Irr(N)}\eta(1)\eta, \lambda]=\lambda(1)=1,$$
which implies that $[\chi_N,\lambda]=e=1$ and $\chi(1)=m=|G:I|$. By Clifford correspondence, there exists $\psi \in \Irr(I)$ such that $\psi^G=\chi$ and $[\psi_N, \lambda] \neq 0$. Since $\chi(1)=|G:I|\psi(1)$, we have $\psi(1)=1$, which implies that $\psi_N=\lambda$. Hence $\chi$ is monomial, as desired.
\end{proof}

\begin{proof}[\bf Proof of Theorem \ref{thma}]
We first show that $(1)$ implies that $(2)$. Let $\chi$ be a $\mathcal{P}$-character of $G$ with respect to a maximal subgroup $M$ and we will show that $\chi$ is monomial. If $\chi$ is principle, then $\chi$ is monomial, as desired. Now we may assume that $\chi$ is non-principle.
	Write $K=\Core_G(M)$. By Lemma~\ref{lem-core}, $K=\Core_G(M) \leq \Ker \chi$. Considering in the quotient group $\overline{G}=G/K$,
	it follows that $\chi \in \Irr(\overline{G})$, moreover, $\chi \in \Irr((1_{\overline{M}})^{\overline{G}})$. Let $\overline{N}=N/K$ be a minimal normal subgroup of $\overline{G}$. The solvability of $G$ implies that $\overline{N}$ is Abelian.  Since $\overline{M}$ is maximal in $\overline{G}$ and core-free, we obtain that $\overline{G}=\overline{N}\overline{M}$ and $\overline{M} \cap \overline{N}=1$. Applying Lemma~\ref{lem-monomial}, $\chi \in \Irr(\overline{G})$ is monomial. Hence $\chi$, as a character of $G$, is also monomial, as desired.

Note that $(2)$ implies $(3)$ trivially by the definition. Now we prove that $(3)$ implies $(1)$. Suppose that the result is false and let 
	$G$ be a counterexample of minimal order. 
	
For any non-trivial normal subgroup $X$ of $G$, we will show that $G/X$ satisfies the hypothesis of Theorem~\ref{thma}. 
	Let $M/X$ be a non-normal maximal subgroup of $G/X$. Then $M$ is also a non-normal maximal subgroup of $G$. By hypothesis, $(1_M)^G$ has a non-principle monomial irreducible constituent $\chi$. By Lemma~\ref{lem-core}, $X \leq \Core_G(M)= \ker \chi$.   
	Hence $\chi \in \Irr(G/X)$ is also a non-principle irreducible constituent of $(1_{M/X})^{G/X}$.  Since $\chi$ is monomial, we may assume 
	$\chi=\eta^G$, where $\eta$ is a linear character of some subgroup $M_1$ of $G$. Then  \[X\leq \text{ker}~\chi= \bigcap\limits_{x\in G}(\text{ker}~ \eta)^x \leq \text{ker}~\eta\leq M_1.\] It follows that $\chi \in {\rm Irr}(G/X)$ is monomial, which is induced by the linear character $\eta \in \Irr(M_1/X)$. The minimality of $G$ implies that $G/X$ is solvable.
	
By following a standard procedure, we can obtain that $\Phi(G)=1$ and $G$ has a unique minimal normal subgroup, $N$ say.
	Clearly $N$ is non-solvable and $N'=N$. Let $H$ be a subgroup of $G$ not containing $N$, chosen to have the largest possible order.
	Then $H \neq 1$ is maximal in $G$ and $G=NH$. Note that $\Core_G(H)=1$ and $H$ is not normal in $G$ otherwise $N \leq H$, contrary to the choice of $H$. By hypothesis, $(1_H)^G$ has a non-principal monomial irreducible constituent, $\psi$ say. By Lemma~\ref{lem-core},  $\text{ker}~\psi=\Core_G(H)=1$.

	Since $\psi$ is monomial, we may assume that $\psi=\mu^G$, where $\mu$ is a linear character of some subgroup $H_1$ of $G$. Then $\psi(1)=|G:H_1|\mu(1)=|G:H_1|$. Since
	\[(1_H)^G(1)=|G:H|\geq 1+\psi(1)= 1+|G:H_1|,\]
	we have that $|H_1|>|H|$. By the maximality of $H$, we see that $N\leq H_1$. As $\mu$ is linear, $H_1' \leq \text{ker}~\mu$. 
	It follows that $1\not=N=N'\leq H_1'\leq\text{ker}~\mu$, and 
	$$N\leq \bigcap\limits_{x\in G}(\text{ker}~ \mu)^x=\text{ker}~\mu^G=\text{ker}~\psi,$$ 
	which contradicts that $\text{ker}~\psi=1$. The proof is complete.
\end{proof}


\begin{proof}[\bf Proof of Theorem \ref{thmgl}]
If $G$ is nilpotent, then it is clear that $\chi(1)^2$ divides $|G : \operatorname{ker}\chi|$ for all $\chi \in \operatorname{Irr}(G)$, and also for all ${\mathcal{P}}$-characters.

Conversely, assume that the result is not true and let $G$ be  a counterexample with minimal order. Then
$\chi(1)^2$ divides $|G : \operatorname{ker}\chi|$ for every ${\mathcal{P}}$-character $\chi \in \operatorname{Irr}(G)$ and $G$ is not nilpotent. 

For each $1 \neq X \unlhd G$ and  for each ${\mathcal{P}}$-character $\chi \in \operatorname{Irr}(G/X)$, $\chi$ is also a $\mathcal{P}$-character
of $G$ with $X \Ker \chi$. By hypothesis, $\chi(1)^2$  divides $|G : \operatorname{ker}\chi|=|G/X: \operatorname{ker}\chi/N|$. The minimality of $G$ implies that $G/X$ is nilpotent. 

Since $G$ is a minimal counterexample, we may assume that $\Phi(G)=1$ and $G$ has a unique minimal normal subgroup $N$. The solvability of $G$ implies that $N$ is an Abelian $p$-subgroup for some prime $p$.
 Choose a maximal subgroup $H$ of $G$ such that $N \not\le H$. Then $G = NH$. Since $N$ is abelian, we have that $N \cap H$ is normal in $N$ as well as in $G$. By the minimality of $N$, it follows that $N \cap H = 1$.

Considering the action of $H$ on $N$ by conjugation,  as $N$ is a $p$-group, $1\ neq \Cent_N(\Oo_p(H))$ is normal in $G$. By the minimality of $N$, $ \Cent_N(\Oo_p(H))=N$, that is, $\Oo_p(H)$ centralizes $N$. The uniqueness of $N$ implies that $\Oo_p(H) \leq \Cent_H(N) = 1$.
Moreover,  we also see that $H$ acts irreducibly and faithfully on $N$. By above, $H \cong G/N$ is nilpotent. Then $\Oo_p(H)=1$ implies that $H$ is a $p'$-group. 

Now, $H$ acts on the dual group $\operatorname{Irr}(N)$ by $h : \lambda \longrightarrow \lambda^h$ for $h \in H$, $\lambda \in \operatorname{Irr}(N)$. By \cite[Theorem B]{Isaacs99}, we may choose $\lambda \in \operatorname{Irr}(N)$ such that $$|\Cent_H(\lambda)| \leq (|H|/p)^{1/p} \leq (|H|/2)^{1/2} < |H|^{1/2},$$ and so the $H$-orbit of $\lambda$ has size $> |H|^{1/2}$.
Note that 
\[
((1_H)^G)_N = ((1_H)_{H\cap N})^N = \sum_{\eta\in\mathrm{Irr}(N)} \eta(1)\eta = \sum_{\eta\in\mathrm{Irr}(N)} \eta.
\]
Hence there exists $\psi \in \Irr((1_H)^G)$ such that \( \lambda \in \Irr(\psi_N)\). 
It follows from Lemma~\ref{lem-monomial} that $\psi(1)=|G:\Cent_G(\lambda)|=|H:\Cent_H(\lambda)|>|H|^{1/2}$, that is, $\psi(1)^2 > |H|$. However, since $N$ is an abelian normal subgroup, by Ito's Theoren (\cite[Theorem 6.15]{Isaacs1}), $\psi(1)$ divides $|G : N| = |H|$, which is coprime to $|N|$.  By hypothesis, $\psi(1)^2$ divides $|G|$. It implies $\psi(1)^2$ divides $|H|$, a contradiction.
\end{proof}

\begin{proof}[\bf Proof of Theorem \ref{thmb}]
 Suppose that $H$ is a normal $\pi$-complement of $G$, and let $\chi$ be a ${\mathcal{P}_\pi}$-character of $G$ with respect to some maximal subgroup $M$ of $G$ such that $|G:M|$ is a $\pi$-number. Since $H$ is a normal Hall $\pi'$-subgroup of $G$, we have that $H\leq M$. By Lemma~\ref{lem-core}, \(H\leq \Core_G(M) \leq \text{ker}~\chi\) and $\chi(1)$ is a  $\pi$-number, as desired.

Conversely, suppose that all ${\mathcal{P}_\pi}$-character degrees of $G$ are 
 $\pi$-number. Let $N$ be a minimal normal subgroup of $G$. By induction hypothesis, $G/N$ has a normal $\pi$-complement $NK/N$, where $K$ is a Hall $\pi'$-subgroup of $G$. Clearly, $N$ is the unique minimal subgroup of $G$. 
 
 Note that $NK \unlhd G$. In what follows, it is sufficient to show that $NK$ has a normal Hall $\pi'$-subgroup. Since \( G \) is \( \pi \)-solvable, we have that \( N \) is a \( \pi' \)-group or elementary abelian \( p \)-group for some prime  \( p \in \pi \). If \( N \) is a \( \pi' \)-group, then \( K=NK \) is a normal Hall \( \pi' \)-subgroup,  as desired.  

Now suppose \( N \) is an elementary abelian \( p \)-group, where \( p \in \pi \).  As $NK \unlhd G$ and $K$ is a Hall $\pi'$-subgroup of $NK$,
by the Frattini argument, \( G = NK\N_G(K) = N\N_G(K) \). If $N \unlhd \N_G(K)$, then \( G = \N_G(K) \), so \( K \unlhd G \), as desired.
Hence we may assume that $N \leq \N_G(K)$.  Write $M=\N_G(K)$ and $G=NM$. As $N \nleq M$, $N \cap M<N$. Since $N$ is Abelian, $N \cap M \unlhd NM=G$. The minimality of $N$ implies that
$N \cap M=1$. Note that
\[
((1_M)^G)_N = ((1_M)_{M\cap N})^N = \sum_{\eta\in\mathrm{Irr}(N)} \eta(1)\eta = \sum_{\eta\in\mathrm{Irr}(N)} \eta.
\]
For each $\lambda \in \Irr(N)$,  there exists an irreducible constituent $\chi$ of \( (1_M)^G \) such that \( \lambda \in \Irr(\chi_N) \). Clearly, \( \chi \) is a  ${\mathcal{P}_\pi}$-character as $|G:M|=|N|$ is a $\pi$-number. By hypothesis, \( \chi(1) \) is a $\pi$-number. 
It follows from Lemma~\ref{lem-monomial} that \( \chi(1) = |G : \Cent_G(\lambda)|=|NM: N\Cent_M(\lambda)|=|M: \Cent_M(\lambda)| \) is a $\pi$-number, and 
$\chi=\psi^G$ for some $\psi \in \Irr(\Cent_G(\lambda))$ with $\psi_N=\lambda$. Hence we see that \( \lambda \) extends to \(\Cent_G(\lambda)=N\Cent_M(\lambda)\).  Since $K$ is a normal Hall $p'$-subgroup of $M$ and $|M: \Cent_M(\lambda)|$ is a $\pi$-number, we obtain that $K \leq \Cent_M(\lambda)$. Hence, \( \lambda \) extends to \(NK\). By \cite[Theorem 3.17]{Isaacs2}, \( NK \) has a normal Hall \( \pi' \)-subgroup $K$, as required.
\end{proof}

\begin{proof}[\textbf{\emph{Proof of Theorem \ref{thmc}}}]
If $G$ is a $\pi$-group, then clearly $\Irr_{\mathcal{P}_{\pi}}(G)=\Irr_{\mathcal{P}}(G)$ holds. Hence it suffices to prove the converse.
Assume that the result is false and let $G$ be a counterexample with minimal order. Then $\Irr_{\mathcal{P}_{\pi}}(G)=\Irr_{\mathcal{P}}(G)$ and $G$ is not a $\pi$-group.

For each $1 \neq X \unlhd G$, we will show that $\Irr_{\mathcal{P}_{\pi}}(G/X)=\Irr_{\mathcal{P}}(G/X)$. In fact, for each $\chi \in \Irr_{\mathcal{P}}(G/X)$, $\chi$ can be viewed as a $\mathcal{P}$-character of $G$ with $X \subseteq \Ker \chi$.  As $\Irr_{\mathcal{P}_{\pi}}(G)=\Irr_{\mathcal{P}}(G)$,  we may assume that $\chi \in \Irr((1_M)^G)$ for some maximal subgroup $M$ of $G$ with $|G:M|$ $\pi$-number. 
By Lemma~\ref{lem-core}, $X \leq \Ker \chi=\Core_G(M) \leq M$. Hence $\chi \in \Irr((1_{M/X})^{G/X})$, which implies that $\chi \in \Irr_{\mathcal{P}_{\pi}}(G/X)$, as desired.  By the minimality of $G$, $G/X$ is a $\pi$-group.

By a standard procedure, $\Phi(G)=1, \Oo_{\pi}(G)=1$ and $G$ has a unique minimal normal subgroup $N$. It follows from Lemma~\ref{lem-core} that
 $$1=\Phi(G)= \bigcap_{\chi \in \Irr_{\mathcal{P}}(G)} \ker \chi.$$
The uniqueness of $N$ implies that there exists $\psi \in \Irr_{\mathcal{P}}(G)$ with $\Ker \psi=1$. By hypothesis, we may assume that $\psi \in \Irr((1_H)^G)$ for some maximal subgroup $H$ of $G$ with $|G:H|$ $\pi$-number.  Since $\Oo_{\pi}(G)=1$ and $G$ is $\pi$-separable, we have $\Oo_{\pi'}(G) \neq 1$. The uniqueness of $N$ implies that $N \leq \Oo_{\pi'}(G)$. Since $|G:H|$ is a $\pi$-number, $N \leq \Oo_{\pi'}(G) \leq H$. It follows from Lemma~\ref{lem-core} that
$1\neq N \leq \Core_G(H)=\Ker \psi=1$, which is a contradiction. The proof is complete.
\end{proof}





\noindent{\bf Acknowledgements}\\

The authors would like to thank Prof. Gabriel Navarro for several useful discussions and comments on this paper. The first author is supported by Guangxi Natural Science Foundation Program (2024GXNSFAA010514). The second author is supported by  Natural Science Foundation of Shanghai 
(24ZR1422800) and National Natural Science Foundation of China (12471018). Part of this work was completed during two authors' visit to the Shenzhen International Center for Mathematics, Southern University of Science and Technology. The authors both express gratitude to Prof. Z. Feng for his hospitality and the generous tea breaks.

\bibliographystyle{amsplain}

\end{document}